\documentclass[11pt,leqno,fleqn]{article}
\usepackage{amsmath}
\usepackage{amsfonts}
\usepackage{psfrag}
\newcommand{\bgcirc}{\raisebox{1.4pt}{\mbox{\scriptsize $\mathbf{\bigcirc}$}}}
\newcommand{\oZ}{{\mathbb{Z}}}
\newcounter{sectie}
\newcommand{\sect}[2]{\refstepcounter{sectie}
\section*{\boldmath \thesectie. #2}%
\label{#1}}
\newcommand{\GG}{{\cal G}}
\newcommand{\oR}{{\mathbb{R}}}
\newcommand{\dy}[2]{%
\refstepcounter{equation}%
\label{#1}%
\begin{list}{}{
\topsep 3mm
\leftmargin 18mm
\rightmargin 0cm
\itemsep 0mm
\listparindent 0mm
\parsep 0mm
\itemsep 0mm
\labelsep 0mm
\labelwidth 18mm
}%
\item[\rm (\theequation)\hfill]
#2
\end{list}%
}
\newcommand{\dyy}[2]{\dy{#1}{\raggedright$\dps#2$}}
\newcommand{\dps}{\displaystyle}
\newcommand{\join}[1]{\raisebox{-.05\height}{\mbox{\hspace*{2pt}\footnotesize$\overset{\hspace*{0.5pt}{#1}}{{\scriptsize\vee}}$\hspace*{2pt}}}}
\newcommand{\rf}[1]{{\rm (\ref{#1})}}
\newcommand{\dez}[1]{\dyz{\raggedright$\displaystyle #1 $}}
\newcommand{\dyz}[1]{%
\refstepcounter{equation}%
\begin{list}{}{
\topsep 3mm
\leftmargin 18mm
\rightmargin 0cm
\itemsep 0mm
\listparindent 0mm
\parsep 0mm
\itemsep 0mm
\labelsep 0mm
\labelwidth 18mm
}%
\item[\rm (\theequation)\hfill]
#1
\end{list}%
}
\newcommand{\thmnn}[1]{\vspace{4mm}\noindent{\bf Theorem.}{\it #1}}
\newcommand{\cornn}[1]{\vspace{4mm}\noindent{\bf Corollary.}{\it #1}}
\newcommand{\dyyz}[1]{\dyz{\raggedright$\dps#1$}}
\newcommand{\thetagraaf}{\raisebox{-.25\height}{\rotatebox{90}{\scalebox{0.02}{\includegraphics{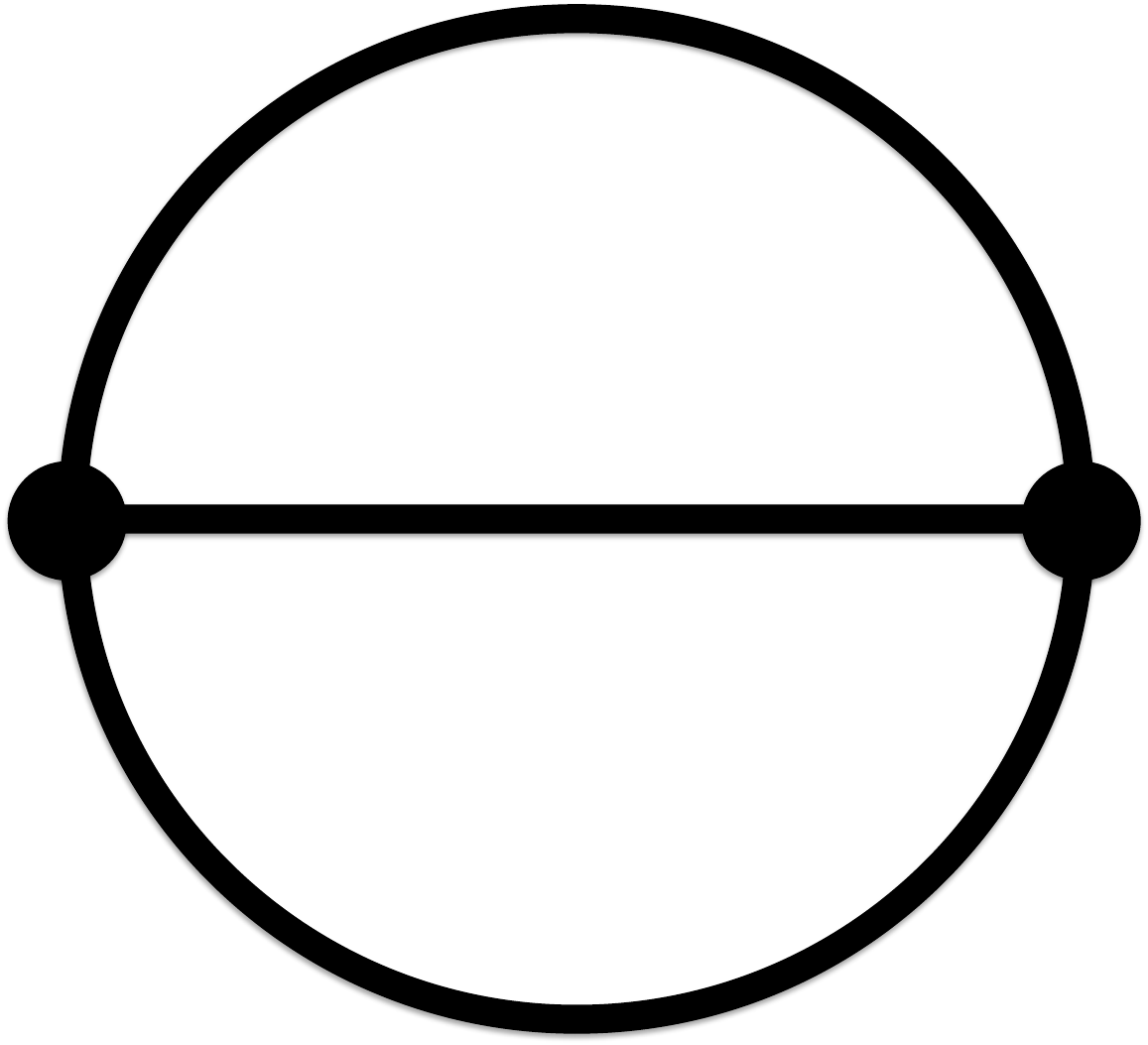}}}}}
\newcommand{\thetagraafw}{\raisebox{-.25\height}{\scalebox{0.02}{\includegraphics{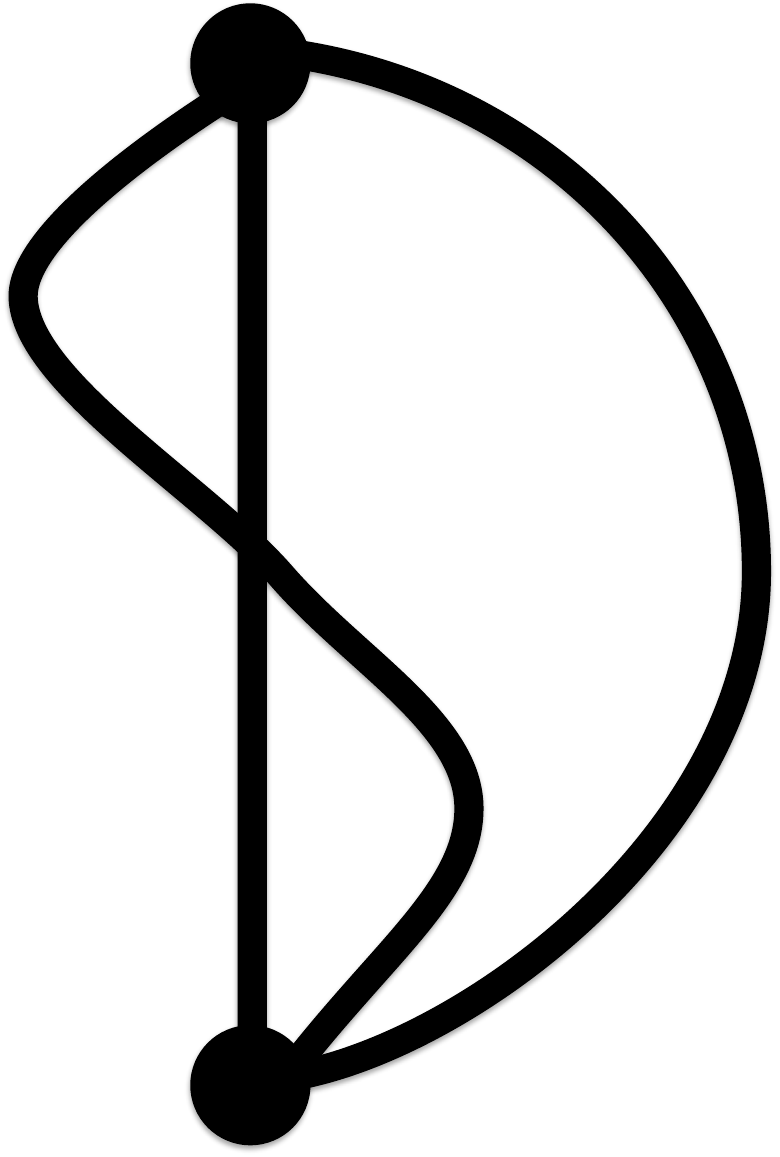}}}}
\newcommand{\pf}{\vspace{3mm}\noindent{\bf Proof.}\ }
\newcommand{\bx}{\hspace*{\fill} \hbox{\hskip 1pt \vrule width 4pt height 8pt depth 1.5pt \hskip 1pt}

\addvspace{4mm}}
\newcounter{hulpstelling}
\newcommand{\lemma}[2]{\refstepcounter{hulpstelling}\vspace{4mm}\noindent{\bf Lemma \thehulpstelling.}\label{#1}{\it #2}}
\newcommand{\de}[2]{\dy{#1}{\raggedright$\displaystyle #2 $}}
\newcommand{\kies}[2]{\mbox{${{#1}\choose{#2}}$}}
\newcommand{\sectz}[1]{\refstepcounter{sectie}
\section*{\boldmath \thesectie. #1}%
}
\newcommand{\kfrac}[2]{\mbox{$\frac{#1}{#2}$}}
\newcommand{\MM}{{\cal M}}
\newcommand{\T}{^{\sf T}}
\newcommand{\minus}{{\hspace*{-1pt}-\hspace*{-1pt}}}
\renewcommand{\phi}{\varphi}
\newcommand{\sgn}{\text{\rm sgn}}
\newcommand{\id}{\text{\rm id}}
\newcommand{\OO}{{\cal O}}
\begin{document}

\begin{center}
{\large\bf ON PARTITION FUNCTIONS FOR 3-GRAPHS

}
\vspace{4mm}
Guus Regts\footnote{ University of Amsterdam.
The research leading to these results has received funding from the European Research Council
under the European Union's Seventh Framework Programme (FP7/2007-2013) / ERC grant agreement
n$\mbox{}^{\circ}$ 339109.},
Alexander Schrijver$\mbox{}^1$,
Bart Sevenster$\mbox{}^1$

\end{center}

\noindent
{\small{\bf Abstract.}
A {\em cyclic graph} is a graph with at each vertex a cyclic order of the edges
incident with it specified.
We characterize which real-valued functions on the collection of cubic cyclic graphs are
partition functions of a real vertex model
(P. de la Harpe, V.F.R. Jones, Graph invariants related to statistical mechanical models:
examples and problems, Journal of Combinatorial Theory, Series B 57 (1993) 207--227).
They are characterized by `weak reflection positivity', which amounts to the positive semidefiniteness of
matrices based on the `$k$-join' of cubic cyclic graphs (for all $k\in\oZ_+$).

Basic tools are the representation theory of the symmetric group and geometric invariant theory,
in particular the Hanlon-Wales theorem on the decomposition of Brauer algebras
and the Procesi-Schwarz theorem on inequalities defining orbit spaces.

}

\sect{22fe15g}{Introduction}

In this paper, a {\em cyclic graph} is an undirected graph where each vertex is equipped with a cyclic order
of the edges incident with it.
A {\em $3$-graph} is a connected cubic cyclic graph.
Loops and multiple edges are allowed.
Also the `vertexless loop' $\bgcirc$ serves as a 3-graph, and may occur as component of a 3-graph.
By $\GG$ and $\GG'$ we denote the collection of 3-graphs and the collection of cubic cyclic graphs, respectively.
So $\GG'$ is the collection of disjoint unions of graphs in $\GG$.

Having fixed this terminology for the current paper, let us stress that the above types of graphs show
up under several names and in several contexts.
Bollob\'as and Riordan [2] refer to Dennis Sullivan for coining the term `cyclic graph'.
Such graphs are also called `graphs with a rotation system'
(cf.\ Gross and Tucker [9]),
`ribbon graphs' (Reshetikhin and Turaev [17]),
or
`fatgraphs' (Milgram and Penner [14]),
and are in one-to-one correspondence with graphs cellularly embedded on a compact oriented surface
(Heffter [12], cf.\ [9]).
Then, by surface duality, cubic cyclic graphs are in one-to-one correspondence with triangulations
of compact oriented surfaces.
The term `3-graph' for a connected cubic cyclic graph was introduced by
Duzhin, Kaishev, and Chmutov [5] (cf.\ [3]), and plays a role, as a variant of chord diagrams and
Jacobi diagrams, in the study of Vassiliev's knot invariants [24].

There is an abundance of literature on invariants for such graphs, introduced to study basic problems
in combinatorics, topology, and theoretical physics.
An important type of invariant is the partition function,
with such basic examples as in the Ising-Potts-model, the Tutte polynomial [23]
(defined for cyclic graphs in [2]),
the Jones polynomial (cf.\ [11]), R-matrices, and Lie algebra weight systems for chord diagrams.

In this paper we focus on partition functions for 3-graphs.
For $n\in\oZ_+$, let $c=(c_{ijk})_{i,j,k=1}^n$ be an element of $((\oR^n)^{\otimes 3})^{C_3}$, which denotes as usual the
linear space of tensors in $(\oR^n)^{\otimes 3}$ that are invariant under the natural action of the
cyclic group $C_3$ on $(\oR^n)^{\otimes 3}$.
Then for any 3-graph $G$ define
\dyy{28se14a}{
f_c(G):=\sum_{\varphi:E(G)\to[n]}\prod_{v\in V(G)}c_{\varphi(e_1)\varphi(e_2)\varphi(e_3)},
}
where, once $v\in V(G)$ is chosen,
$e_1,e_2,e_3$ denote the edges incident with $v$, in cyclic order.
(As usual, $[n]:=\{1,\ldots,n\}$.)
Following the terminology of de la Harpe and Jones [11], $f_c$ is the {\em partition function} of the {\em vertex model}
$c=(c_{ijk})_{i,j,k=1}^n$.
Alternatively, a vertex model is called an {\em edge-coloring model} ([22]).
Note that $f_c(\bgcirc)=n$.

In this paper we will characterize for which functions $f$ on $\GG$ there exist $n\in\oZ_+$
and $c\in((\oR^n)^{\otimes 3})^{C_3}$ with $f=f_c$.
They are exactly characterized by a form of `reflection positivity', a notion which
roots in quantum field theory (cf.\ [7]) and has been defined in various ways for graph parameters
(see e.g. [6],[22],[20]).
Here we adopt a weaker version of it based on the $k$-join of cubic cyclic graphs,
which is a restricted variant of the glueing operation of graphs with open ends considered by
Szegedy [22].
It yields a weaker condition, and therefore it makes the characterization stronger.

To be precise, let $\oR[\GG]$ denote the commutative $\oR$-algebra freely generated by the collection
$\GG$ of 3-graphs.
Any function from $\GG$ to any $\oR$-algebra can be extended uniquely to an algebra homomorphism on
$\oR[\GG]$.
We identify the product $G_1\cdots G_k$ of $3$-graphs in $\oR[\GG]$ with the disjoint union
of $G_1,\ldots,G_k$, which is a cubic cyclic graph.
So the collection $\GG'$ of cubic cyclic graphs corresponds to the set of monomials in $\oR[\GG]$.

For $G$ and $H$ in $\GG'$, the {\em $k$-join} $G\join{k}H$ is the element of $\oR[\GG]$ obtained
as follows.
Consider the disjoint union of $G$ and $H$.
Choose distinct vertices $u_1,\ldots,u_k$ of $G$ and distinct
vertices $v_1,\ldots,v_k$ of $H$, and, for each $i=1,\ldots,k$,
\dy{28se14b}{
replace~~~~%
$\raisebox{-.3\height}{\scalebox{0.10}{\includegraphics{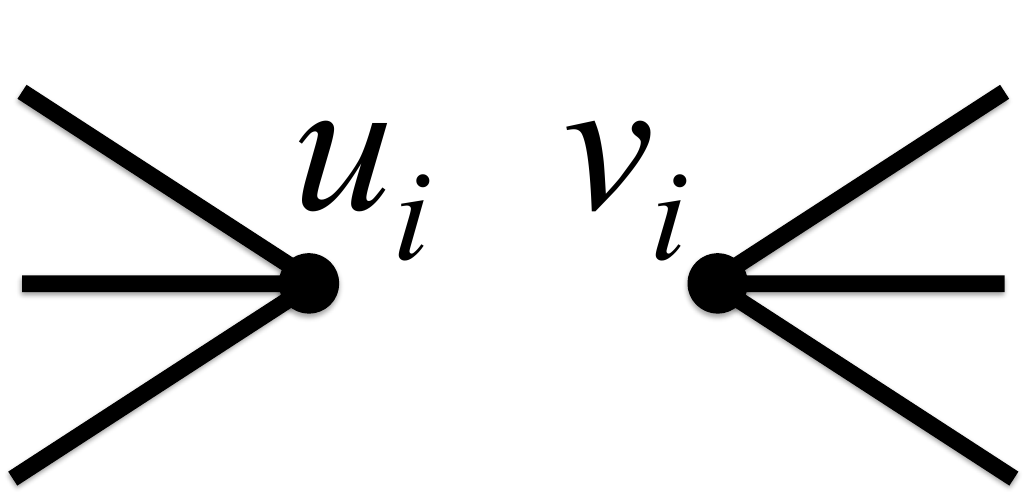}}}$~~~~by%
~~~~$\frac{1}{3}\big(
~~
\raisebox{-.4\height}{\scalebox{0.10}{\includegraphics{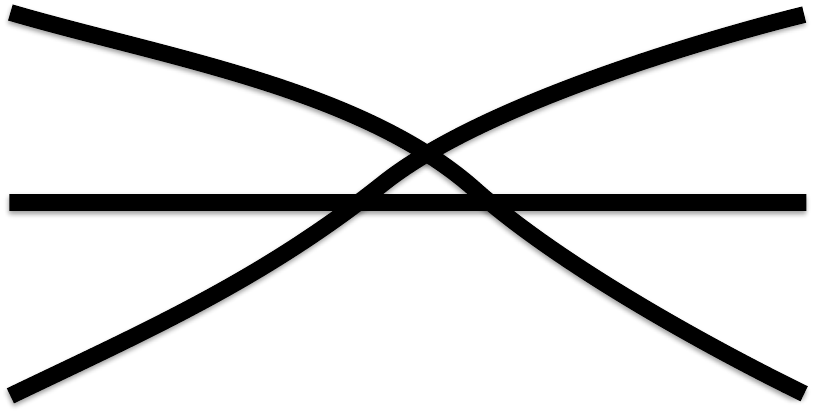}}}
~~+~~
\raisebox{-.4\height}{\scalebox{0.10}{\includegraphics{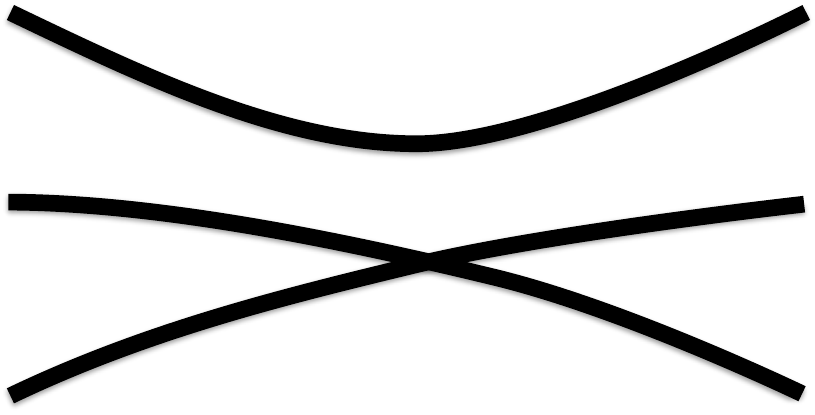}}}
~~+~~
\raisebox{-.4\height}{\scalebox{0.10}{\includegraphics{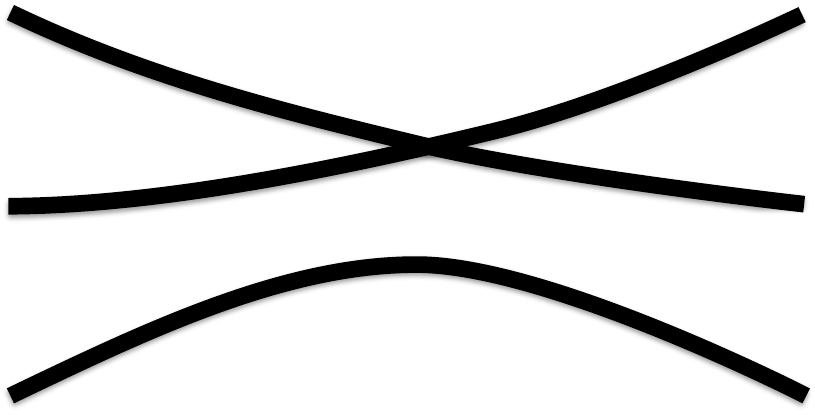}}}
~~
\big)$~.
}
Here and below, the cyclic order at a vertex is that given by the clockwise orientation.
(Thus the new connections in \rf{28se14b} obey the cyclic orders at $u_i$ and $v_i$: choosing total orders compatible
with the cyclic orders,
the first edge at $u_i$ is connected with the first edge at $v_i$,
the second edge at $u_i$ is connected with the second edge at $v_i$,
and the third edge at $u_i$ is connected with the third edge at $v_i$.
There are three different ways of doing so.)
Then $G\join{k}H$ is obtained by adding up these elements of $\oR[\GG]$ over all choices of $u_1,\ldots,u_k$
and $v_1,\ldots,v_k$.
(We add up over all possible orderings of the vertices $u_1,\ldots,u_k$ and $v_1,\ldots,v_k$.)
(The `$k$-join' could be described dually in terms of
surgery of triangulated surfaces.
However, we will not pursue this visualisation, as the $k$-join seems easier to handle
in the cyclic graph setting.)

Then $f:\GG\to\oR$ is called {\em weakly reflection positive} if, for each $k$, the $\GG'\times\GG'$ matrix
\dez{
M_{f,k}:=(f(G\join{k}H))_{G,H\in\GG'}
}
is positive semidefinite
(that is, each finite principal submatrix is positive semidefinite).

We can extend $G\join{k}H$ bilinearly to a bilinear function $\oR[\GG]\times\oR[\GG]\to\oR[\GG]$.
Then weak reflection positivity means that
$f(\gamma\join{k}\gamma)\geq 0$ for each $\gamma\in\oR[\GG]$ and each $k\in\oZ_+$.

\thmnn{
Any function $f:\GG\to\oR$ is the partition function of some real vertex model
if and only if $f$ is weakly reflection positive.
}

\medskip
We prove this theorem in Section \ref{17fe15a}.
It is not hard to show that if $f$ is a partition function, then $f=f_c$ for some unique
$n\in\oZ_+$ and $c\in((\oR^n)^{\otimes 3})^{C_3}$, up to the natural action of the real orthogonal
group $O(n)$ on $c$ (which leaves $f_c$ invariant) --- see Section \ref{22fe15b}.

We derive a direct consequence of the theorem that considers {\em weight systems}.
They play a key role in the study of Vassiliev's invariants
for classifying
the finite-type invariants for knots of Vassiliev [24] (through the Kontsevich integral [13])
and for integral homology 3-spheres (Ohtsuki [15]).

For 3-graphs, a (real-valued) {\em weight system} is a function $f:\GG\to\oR$
which is {\em antisymmetric}\/:
$f(H)=-f(G)$ if $H$ arises from $G$ by reversing the orientation at one of its vertices, and
satisfies the {\em IHX-equation} (which roots in work of Bar-Natan, cf.\ [1]):
\dyyz{
f(~\raisebox{-.2\height}{\scalebox{0.08}{\includegraphics{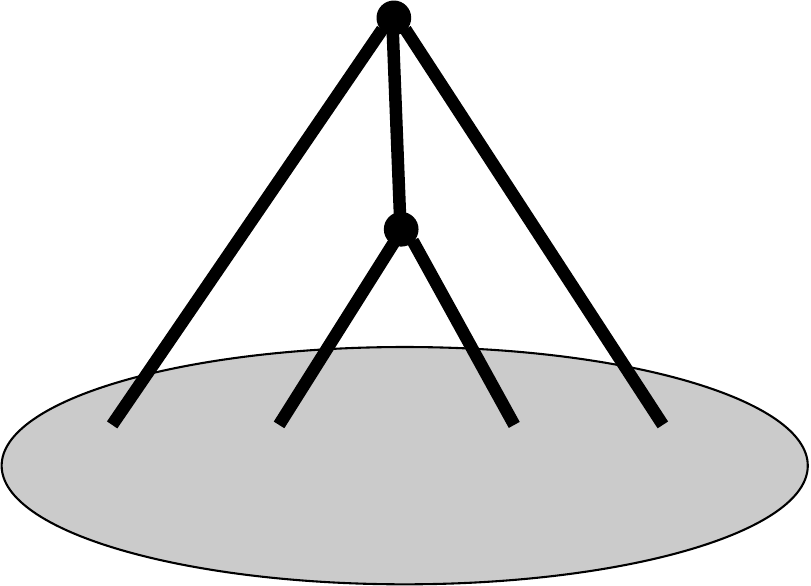}}}~)
~~=~~
f(~\raisebox{-.3\height}{\scalebox{0.08}{\includegraphics{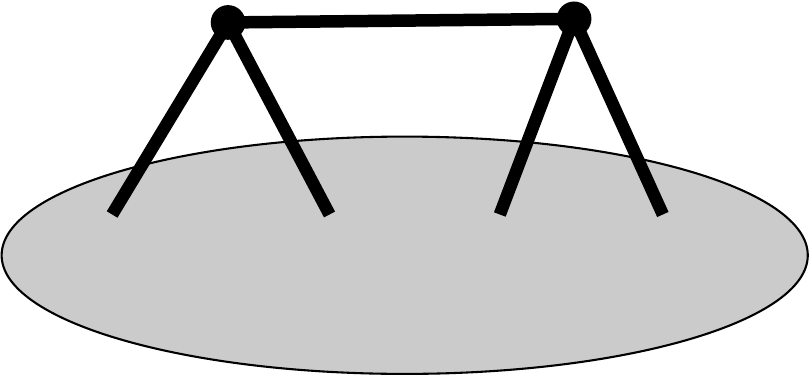}}}~)
~~-~~
f(~\raisebox{-.3\height}{\scalebox{0.08}{\includegraphics{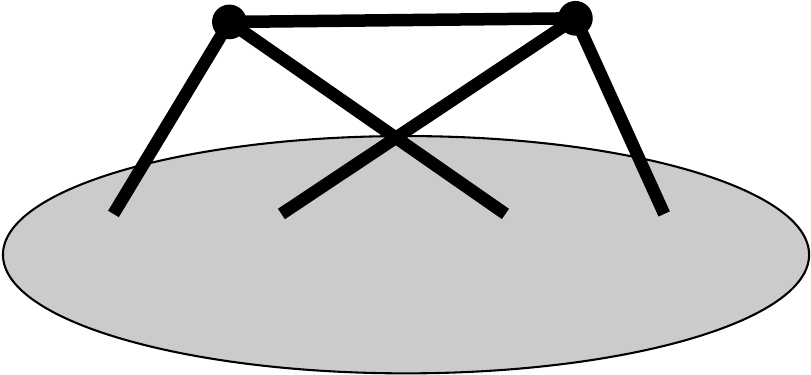}}}~).
}
Key instances of weight systems are the {\em Lie algebra weight systems}\/: the partition functions $f_c$
of the structure tensor $c$ of a finite-dimensional Lie algebra $\mathfrak{g}$, expressed in a basis that
is orthonormal with respect to some symmetric ad-invariant bilinear form on $\mathfrak{g}$.

\cornn{
A function $f:\GG\to\oR$ is a Lie algebra weight system
if and only if $f$ is weakly reflection positive and satisfies
$f(\thetagraafw)=-f(\thetagraaf)$
and
$f(\raisebox{-.25\height}{\scalebox{0.02}{\includegraphics{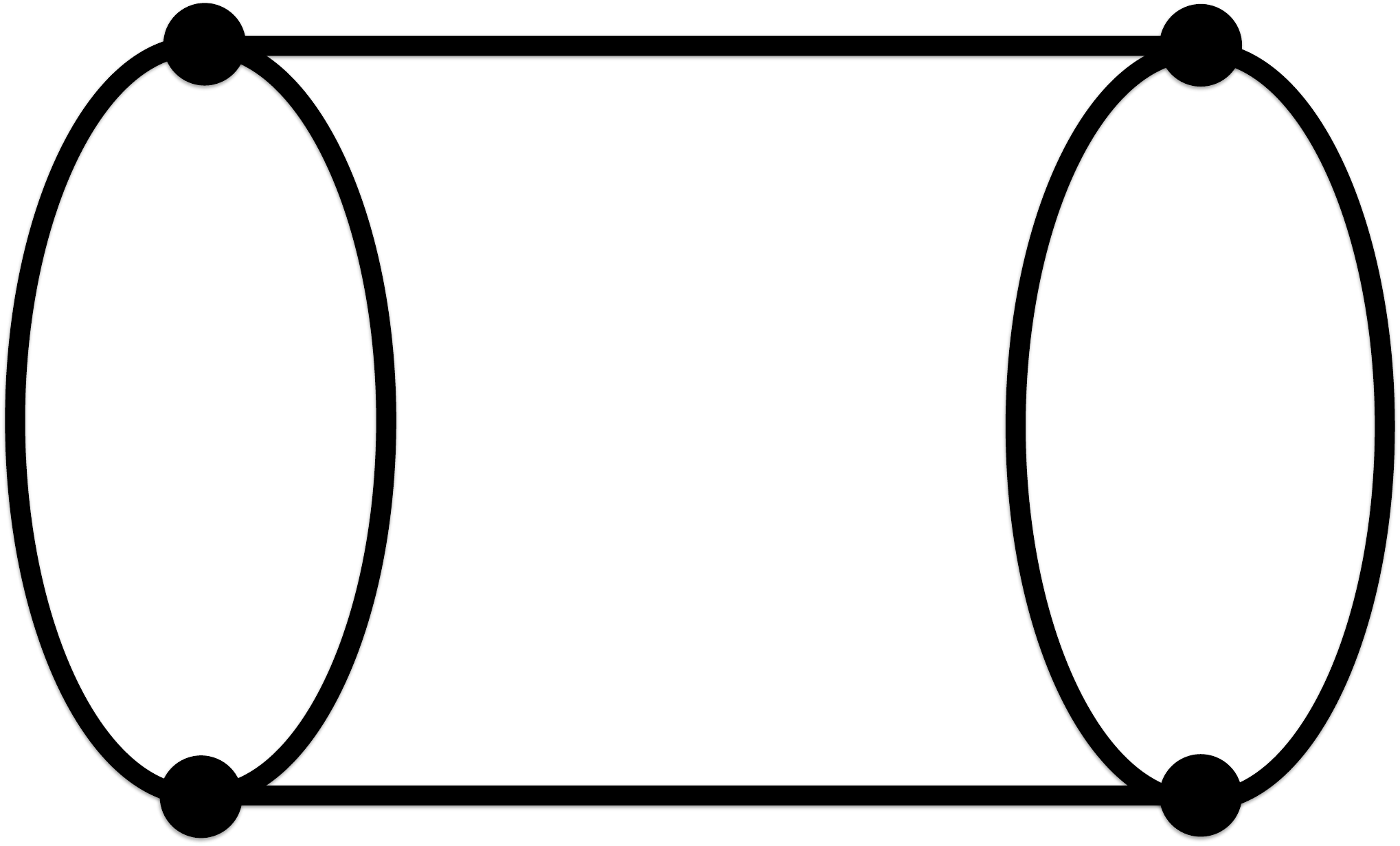}}})=
2f(\raisebox{-.2\height}{\scalebox{0.02}{\includegraphics{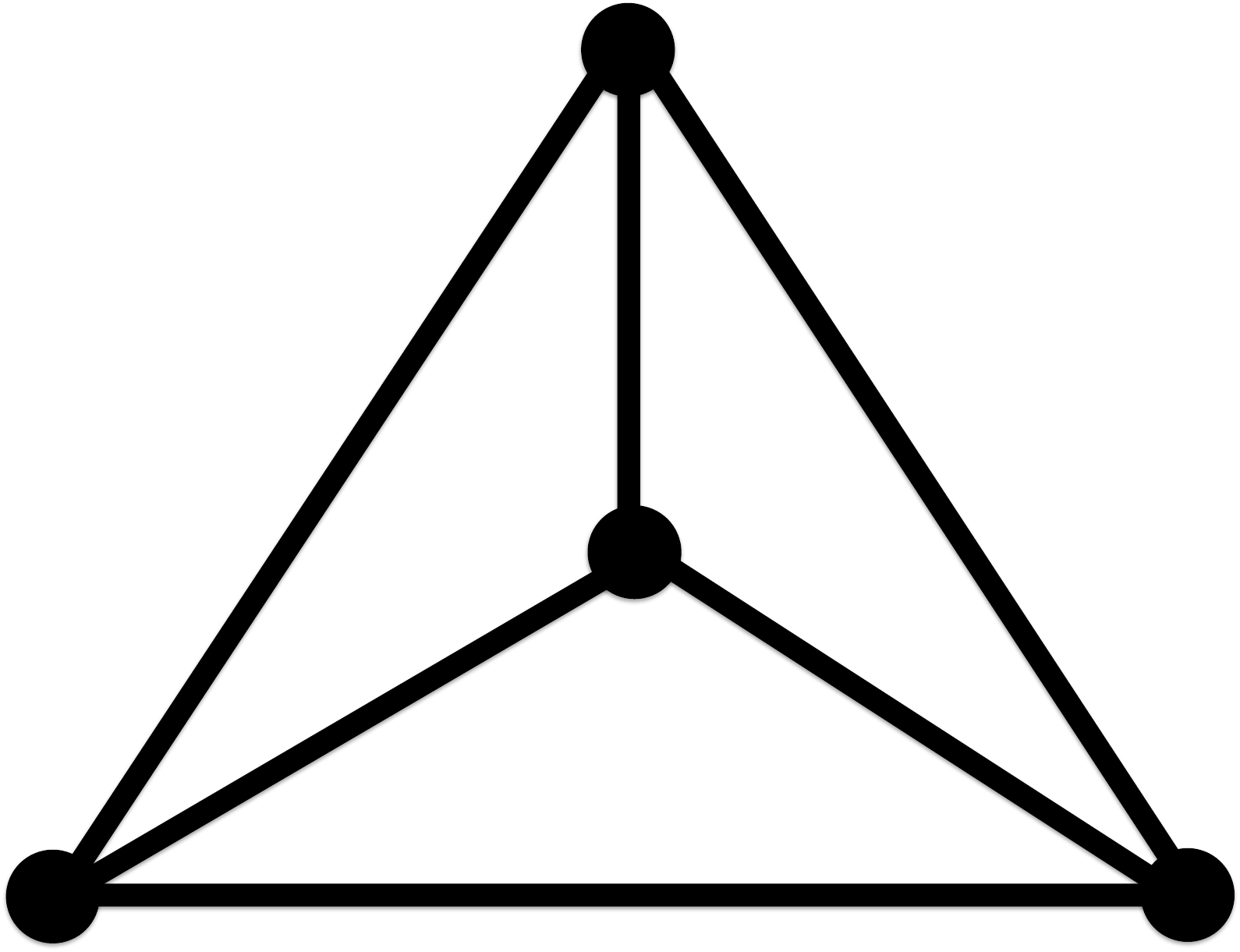}}})$.
}

\pf
This follows from the theorem, as for any $n$ and any $c\in((\oR^n)^{\otimes 3})^{C_3}$, if
$f_c(\thetagraafw)=-f_c(\thetagraaf)$
then $c$ is an alternating tensor, as $f_c(\thetagraafw)=-f_c(\thetagraaf)$
is equivalent to
\dy{29jl16a}{
$\dps\sum_{i,j,k}(c_{ijk}+c_{ikj})^2=0$, and hence to: $c_{ikj}=-c_{ijk}$ for all $i,j,k$.
}
Indeed, 
$f_c(\thetagraafw)=\sum_{i,j,k}c_{ijk}c_{ijk}$, as the cyclic order of the edges at both vertices
of $\thetagraafw$ are equal.
Moreover,
$f_c(\thetagraaf)=\sum_{i,j,k}c_{ijk}c_{kji}$, as the cyclic order of the edges at both vertices
of $\thetagraaf$ are opposite.
Hence (as $c$ is $C_3$-invariant)
\dyyz{
\sum_{i,j,k}(c_{ijk}+c_{ikj})^2=
\sum_{i,j,k}(c_{ijk}^2+2c_{ijk}c_{ikj}+c_{ikj}^2)=
2f_c(\thetagraafw)+2f_c(\thetagraaf)=0.
}
Thus we have \rf{29jl16a}.

Similarly, if $c$ is alternating, then
$f_c(\raisebox{-.25\height}{\scalebox{0.02}{\includegraphics{h4.pdf}}})=
2f_c(\raisebox{-.2\height}{\scalebox{0.02}{\includegraphics{k4.pdf}}})$
is equivalent to $c$ being
the structure tensor of some Lie algebra expressed in a basis that is orthonormal
with respect to some symmetric ad-invariant bilinear form, as it is equivalent to
\dyyz{
\sum_{i,j,k,l}\big(\sum_a(c_{ija}c_{akl}+c_{ila}c_{ajk}+c_{ika}c_{alj})\big)^2=0,
}
and hence to: $\sum_a(c_{ija}c_{akl}+c_{jka}c_{ail}+c_{kia}c_{ajl})=0$ for all $i,j,k,l$.
This last set of equations amounts to the Jacobi identity.
\bx

The theorem is proved using the decomposition of Brauer algebras as given by Hanlon and Wales [10],
the first fundamental theorem of invariant theory, and
the characterization of orbit spaces by inequalities of Procesi and Schwarz [16].

Compared with previous work on this type of issue, the present paper considers $k$-joins
and uses the Procesi-Schwarz theorem, instead of joining graphs with open ends and using
the real Nullstellensatz as in Szegedy [22].
Compared with [19], instead of general graphs the present paper is considering
3-graphs, for which we need to apply deeper representation theory (of the symmetric group)
to derive that $f(\bgcirc)$ is an integer.
Furthermore, a `$k$-join lemma' is given below that simplifies the proof.
The complex case, as studied in [4],[21], demands different
conditions and machinery, and requires (so far) the dimension of the vertex model to be specified in the theorem.

We do not know whether the positive semidefiniteness condition
can be further relaxed to a variant of the $k$-join in which we add in \rf{28se14b} also the
three anti-cyclic connections, each with a minus sign, in line with the vertex product of
$3$-graphs of Duzhin, Kaishev, and Chmutov [5] (cf.\ [3]).
If the integrality of $f(\bgcirc)$ can be derived also for this even weaker form of reflection positivity,
the rest of the proof and hence the theorem will be maintained.

We now first prove three lemmas (in Sections \ref{22fe15e}--\ref{22fe15f}), with which
the proof of the theorem
in Section \ref{17fe15a} follows by a concise series of arguments based on invariant theory.
In Section \ref{22fe15b} we show the uniqueness of the vertex model $c$.

\sect{22fe15e}{A $k$-join lemma}

In the following lemma, $\vartheta$ denotes the 3-graph $\thetagraafw$, and
$\vartheta^i$ is the $i$-th power of $\vartheta$,
that is, the disjoint union of $i$ copies of $\thetagraafw$.

\lemma{13ap08f}{
For any $k$ and any $G\in\GG'$ with $n$ vertices:
\de{22fe15d}{
\kies{n}{k}G
=
2^{-k}k!^{-2}
\sum_{i=0}^k(-1)^{k-i}\kies{k}{i}(G\join{k}\vartheta^i)\vartheta^{k-i}.
}
}

\pf
For each $i$, let $G\underline{\join{k}}\vartheta^i$ be equal to the sum describing $G\join{k}\vartheta^i$
in Section \ref{22fe15g} (with $H:=\vartheta^i$) restricting the summation to those
$v_1,\ldots,v_k$ where each component of $\vartheta^i$ contains at least one vertex among
$v_1,\ldots,v_k$.
So for each $i$,
$G\join{k}\vartheta^i
=
\sum_{j=0}^i
\kies{i}{j}
(G\underline{\join{k}}\vartheta^j)
\vartheta^{i-j}$.
Hence
\dyyz{
\sum_{i=0}^k
(-1)^{k-i}
\kies{k}{i}
(G\join{k}\vartheta^i)\vartheta^{k-i}
=
\sum_{i=0}^k
(-1)^{k-i}
\kies{k}{i}
\sum_{j=0}^i
\kies{i}{j}
(G\underline{\join{k}}\vartheta^j)
\vartheta^{k-j}
=
\sum_{j=0}^k
\kies{k}{j}
(G\underline{\join{k}}\vartheta^j)\vartheta^{k-j}
\sum_{i=j}^k(-1)^{k-i}\kies{k-j}{k-i}
=
G\underline{\join{k}}\vartheta^k
=
2^kk!^2\kies{n}{k}G,
}
the last equality because
$u_1,\ldots,u_k$ can be chosen in $\kies{n}{k}k!$ ways and $v_1,\ldots,v_k$ in $2^kk!$ ways, while each
term of $G\underline{\join{k}}\vartheta^k$ is equal to $G$.
\bx

\sectz{Integrality of $f(\bigcirc)$}


\lemma{22fe15a}{
If $f:\GG\to\oR$ is weakly reflection positive, then $f(\bgcirc)\in\oZ_+$.
}

\pf
Let $f:\GG\to\oR$ be weakly reflection positive.
A direct computation shows
\dyyz{
(\thetagraafw-\thetagraaf)\join{2}(\thetagraafw-\thetagraaf)=\kfrac23\bgcirc(\bgcirc-1)(\bgcirc-2).
}
By the weak reflection positivity of $f$ this implies
$f(\bgcirc)(f(\bgcirc)-1)(f(\bgcirc)-2)\geq 0$, hence $f(\bgcirc)\geq 0$.
To prove that $f(\bgcirc)$ is integer, define $k:=\lceil f(\bgcirc)\rceil+1$.

Let $\MM$ be the set of perfect matchings on $[6k]$. 
To each $M\in \MM$ we can associate a graph $G_M\in\GG'$ on $[2k]$ by identifying the vertices
$3j-2,3j-1,3j$ of $M$ for $j\in[2k]$, with this cyclic order at $j$.
So
\dyz{
$\raisebox{-.3\height}{\scalebox{0.15}{\includegraphics{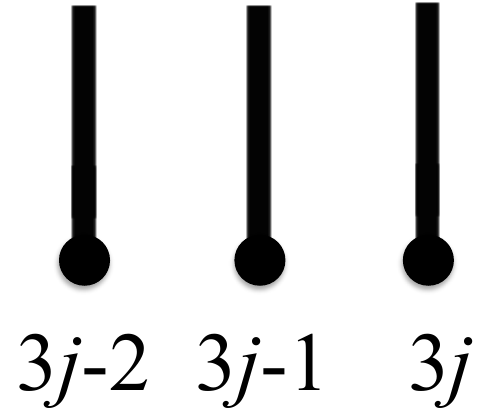}}}$
~~
becomes
~~
$\raisebox{-.3\height}{\scalebox{0.15}{\includegraphics{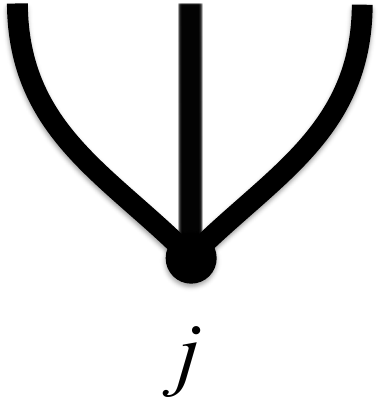}}}$~.
}

For all $M,N\in\MM$, $G_M\join{2k}G_N$ is a polynomial in $\bgcirc$,
since both $G_M$ and $G_N$ have $2k$ vertices.
To describe this polynomial, we consider the natural action of the symmetric group $S_{6k}$ on $\MM$ as:
$\pi\cdot M=\{\pi(e)\mid e\in M\}$ for $M\in\MM$ and $\pi\in S_{6k}$.
This induces an action on $\oR^\MM$ and makes $\oR^\MM$ an $S_{6k}$-module.

For $j\in [2k]$, let $B_j$ be the group of cyclic permutations of $\{3j-2,3j-1,3j\}$,
and define $B:=B_1B_2\cdots B_{2k}$.
Let $D$ be the group of permutations $d\in S_{6k}$ for which there exists $\pi\in S_{2k}$ such that
$d(3j-i)=3\pi(j)-i$ for each $j=1,\ldots,2k$ and $i=0,1,2$.
Set $Q:=BD$, which can be seen to be a group again.

For $M,N\in \MM$, let $c(M,N)$ denote the number of connected components of $M\cup N$. 
Then, by definition of the operation $\join{2k}$, we have
\dyy{eq:vee=union}{
G_M \join{2k} G_N=(2k)!3^{-2k}\sum_{q\in Q}\bgcirc^{c(M,q\cdot N)}.
}

For $\pi\in S_{6k}$, let $P_{\pi}$ be the $\MM\times\MM$ permutation matrix corresponding to $\pi$;
that is, $P_{\pi}w=\pi\cdot w$ for each $w\in\oR^\MM$.
For any $x\in\oR$, let $A(x)$ and $A^Q(x)$ be the $\MM\times\MM$ matrices defined by 
\dyz{
$(A(x))_{M,N} :=x^{c(M,N)}$
~~and~~
$\dps A^Q(x):=\sum_{q\in Q}A(x)P_q$,
}
for $M,N\in \MM$.
Note that each $P_{\pi}$ commutes with $A(x)$, as for all $M,N\in\MM$
one has $c(\pi\cdot M,\pi\cdot N)=c(M,N)$, implying $A(x)=P_{\pi}\T A(x)P_{\pi}=P_{\pi}^{-1}A(x)P_{\pi}$.

Define
\dyyz{
\mu(x):=\prod_{i=0}^{k-1}(x-i)(x-i+2)(x+2i+4).
}
Then we claim that
\dy{prop:eigen}{
$|Q|\mu(x)$ is an eigenvalue of $A^Q(x)$. 
}
This implies the lemma, since by
\eqref{eq:vee=union}, $A^Q(f(\bgcirc))_{M,N}=(2k!)^{-1}3^{2k}f(G_M\join{2k}G_N)$.
Hence, by the weak reflection positivity of $f$, $A^Q(f(\bgcirc))$ is positive semidefinite.
So $\mu(f(\bgcirc))\geq 0$, hence, as $k-1=\lceil f(\bgcirc)\rceil$ and as $k-1$ is the largest
zero of $\mu(x)$, with multiplicity 1, we know $f(\bgcirc)=k-1$.

To prove \rf{prop:eigen}, we will give an eigenvector $u$ of $A^Q(x)$ belonging to $|Q|\mu(x)$.
We derive $u$ from the eigenvector $v$ of $A(x)$ belonging to $\mu(x)$ as described by the
formula for $h_{\lambda}(x)$ in Theorem 3.1 of
Hanlon and Wales [10] as follows, using the representation theory of $S_{6k}$ (cf.\ Sagan [18]).

Consider the following Young tableau, associated to the partition
$(2k+4,\underbrace{4,\ldots,4}_{k-1})$ of $6k$:
\renewcommand{\arraystretch}{1.2}
\dyyz{
T:=\begin{array}{|c|c|c|c|c|c|c|c|c|c|c|c|c|}
\hline
1&\overline 1&2&\overline 2&3&\overline 3&6&\overline 6&9&\overline 9&\cdots&3k&\overline{3k}\\
\hline
4&\overline 4&5&\overline 5\\
\cline{1-4}
7&\overline 7&8&\overline 8\\
\cline{1-4}
\vdots&\vdots&\vdots&\vdots\\
\cline{1-4}
3k\minus 2&\overline{3k\minus 2}&3k\minus 1&\overline{3k\minus 1}\\
\cline{1-4}
\end{array},
}
where $\overline i:=3k+i$ for $i\in[3k]$.

Let $F$ be the perfect matching in $\MM$ with edges $\{i,\overline i\}$, for $i\in[3k]$.
For $i=1,\ldots,4$, let $K_i$ denote the set of elements in the $i$-th column of $T$ and let $C_i$ be
the subgroup of $S_{6k}$ that permutes the elements of $K_i$.
Then $C$ is the group $C_1C_2C_3C_4$.
Similarly, for $i=1,\ldots,k$, let $R_i$ be the subgroup of $S_{6k}$ that permutes the numbers in row $i$ of $T$
and leaves all other numbers fixed,
and $R$ is the group $R_1\cdots R_k$.
Define $v$ and $u$ in $\oR^\MM$ by
\dyz{
$\dps v:=\sum_{c\in C,r\in R} \sgn(c)cr\cdot F$
~~and~~
$\dps u:=\sum_{q\in Q} q\cdot v$,
}
identifying an element of $\MM$ with the corresponding basis vector in $\oR^\MM$.
By Theorem 3.1 of [10] and its proof, $v$ is an eigenvector of $A(x)$ with eigenvalue $\mu(x)$.
Hence
\dyyz{
A^Q(x)u=\sum_{q',q\in Q}AP_{q'}P_{q}v=\sum_{q',q\in Q}P_{q'}P_{q}Av=\mu(x)\sum_{q',q\in Q}P_{q'}P_{q}v=|Q|\mu(x) u.
}
So to prove \rf{prop:eigen}, and hence the lemma, it suffices to show that $u$ is nonzero.
To this end we show that the coefficient $u_F$ of $F$ in $u$ is nonzero.
Note that
\dyyz{
u_F=\sum_{q\in Q}(q\cdot v)_F=\sum_{q\in Q}\sum_{c\in C,r\in R}\sgn(c)(qcr\cdot F)_F
=
\sum_{q\in Q,c\in C,r\in R\atop qcr\cdot F=F}\sgn(c).
}
So it suffices to show that for any $q\in Q$, $c\in C$, and $r\in R$, if $qcr\cdot F=F$ then $\sgn(c)=1$.
As $Q$ is a group, equivalently it suffices to show for any $q\in Q$, $c\in C$, $r\in R$:
\dyz{
if $q\cdot F=cr\cdot F$, then $\sgn(c)=1$.
}
Choose $q\in Q$, $c\in C$, and $r\in R$ with $q\cdot F=cr\cdot F$.
Let $c=c_1c_2c_3c_4$ with $c_i\in C_i$ ($i=1,\ldots,4$) and define $M:=q\cdot F$.
Let $\zeta\in S_{6k}$ be defined by $\zeta(i):=i+1$ if $3$ does not divide $i$ and $\zeta(i):=i-2$ if 3 divides $i$.
So $\zeta^3=\id$ and $\zeta\cdot F=F$.
Moreover, $\zeta q=q\zeta$ (since $\zeta b=b\zeta$ and $\zeta d=d\zeta$ for all $b\in B$ and $d\in D$).
Hence $\zeta\cdot M=M$.

Let $\phi(i):=\overline i$ for $i\in [3k]$.
We show that for each $a\in K_1$:
\dyy{14fe15a}{
c_2\phi c_1^{-1}(a)=\zeta^{-1}c_4\phi c_3^{-1}\zeta(a).
}
This implies $\sgn(c_2c_1^{-1})=\sgn(c_4c_3^{-1})$, and hence $\sgn(c)=1$.

As both $c_2\phi c_1^{-1}$ and $\zeta^{-1}c_4\phi c_3^{-1}\zeta$ are bijections $K_1\to K_2$, it suffices to show
\rf{14fe15a} for all $a\in K_1\setminus\{c_1(1)\}$.
Therefore, choose $a\in K_1$ with $i:=c_1^{-1}(a)\neq 1$.
Let $b:=c_2(\overline i)=c_2\phi c_1^{-1}(a)$.
Note that $b\in K_2$, $\zeta(a)\in K_3$, and $\zeta(b)\in K_4$.
We must show that $c_4^{-1}\zeta(b)=\phi c_3^{-1}\zeta(a)$, that is,
$c_3^{-1}\zeta(a)$ and $c_4^{-1}\zeta(b)$ belong to the same row of $T$.

First assume that $\{i,\overline i\}\in r\cdot F$.
Then $\{a,b\}\in cr\cdot F=M$,
hence, by the $\zeta$-invariance of $M$, $\{\zeta(a),\zeta(b)\}\in M$.
So $\{c_3^{-1}\zeta(a),c_4^{-1}\zeta(b)\}$ belongs to $r\cdot F$, and hence
it is contained in a single row of $T$.

Second assume that $\{i,\overline i\}\not\in r\cdot F$.
Since $i\neq 1$, this implies that $i$ and $\overline i$ are matched in $r\cdot F$
with elements of $K_3\cup K_4$.
So $a$ and $b$ are matched in $M$ with elements of $K_3\cup K_4$.
Hence, by the $\zeta$-invariance of $M$, $\zeta(a)$ and $\zeta(b)$ are matched in $M$ with elements
of $\zeta(K_3\cup K_4)$, which is the first row of $T$ outside $K_1\cup K_2\cup K_3\cup K_4$.
So $c_3^{-1}\zeta(a)$ and $c_4^{-1}\zeta(b)$ are matched in $r\cdot F$ with elements of
the first row of $T$, and hence they both also belong to the first row of $T$.
\bx

\sect{22fe15f}{The polynomial $p_n(G)$}

Choose $n\in\oZ_+$ and let $W$ be the linear space
\dyyz{
W:=((\oR^n)^{\otimes 3})^{C_3}.
}
As usual, $\OO(W)$ denotes the algebra of polynomials on $W$.
For each $3$-graph $G$, define the polynomial $p_n(G)\in\OO(W)$ by
$p_n(G)(c):=f_c(G)$ for any $c\in W$ (defined in \rf{28se14a}).
This can be extended uniquely to an algebra homomorphism $p_n:\oR[\GG]\to\OO(W)$.

For any $q\in\OO(W)$, let $dq$ be its derivative, being an element of $\OO(W)\otimes W^*$.
So $d^kq\in\OO(W)\otimes (W^*)^{\otimes k}$.
Note that the standard inner product on $\oR^n$ induces an inner product on $W$, hence on $W^*$, and
hence it induces a product
$\langle.,.\rangle:
(\OO(W)\otimes (W^*)^{\otimes k})\times
(\OO(W)\otimes (W^*)^{\otimes k})\to\OO(W)$,
by
\dyyz{
\langle p\otimes f_1\otimes\cdots\otimes f_k, q\otimes g_1\otimes\cdots\otimes g_k\rangle=
pq\langle f_1\otimes\cdots\otimes f_k,g_1\otimes\cdots\otimes\otimes g_k\rangle
=
pq\langle f_1,g_1\rangle\ldots\langle f_k,g_k\rangle.
}

The following lemma will be used several times in our proof.

\lemma{13ap08d}{
For all $G,H\in\GG'$ and all $k,n\in\oZ_+$:
\de{13ap08c}{
p_n(G\join{k}H)=\langle d^kp_n(G),d^kp_n(H)\rangle.
}
}

\pf
Let $b_1,\ldots,b_n$ be the standard basis of $\oR^n$, with dual basis $b^*_1,\ldots,b^*_n$.
For $i,j,k=1,\ldots,n$, let $y_{ijk}$ be the element $b_i^*\otimes b_j^*\otimes b_k^*|_W$ of $W^*$.

Consider some $G\in\GG'$.
For $\varphi:E(G)\to[n]$ and $v\in V(G)$, denote
\dez{
\widehat\varphi_v:=y_{\varphi(e_1)\varphi(e_2)\varphi(e_3)},
}
where $e_1,e_2,e_3$ are the edges incident with $v$, in order.
Then
\dez{
p_n(G)=\sum_{\varphi:E(G)\to[n]}\prod_{v\in V(G)}\widehat\varphi_v.
}
Hence $d^kp_n(G)$ expands as:
\dyyz{
d^kp_n(G)=
\sum_{\varphi:E(G)\to [n]}
\sum_{u_1,\ldots,u_k\in V(G)}
\big(
\prod_{v\in V(G)\setminus\{u_1,\ldots,u_k\}}
\widehat\varphi_v
\big)
\otimes
\widehat\varphi_{u_1}
\otimes\cdots\otimes
\widehat\varphi_{u_k},
}
with $u_1,\ldots,u_k$ taken distinct.
Now we claim that for all functions $i,j:[3]\to[n]$,
\de{14fe15b}{
\langle y_{i(1)i(2)i(3)},y_{j(1)j(2)j(3)}\rangle
=
\kfrac13|\{\pi\in C_3\mid
j(s)=i(\pi(s))\text{ for }s\in[3]\}|.
}
Indeed, for each $i:[3]\to[n]$ and $x\in W$, by the $C_3$-invariance of $x$:
\dez{
y_{i(1)i(2)i(3)}(x)=
\langle b_{i(1)}\otimes b_{i(2)}\otimes b_{i(3)},x\rangle
=
\langle\kfrac13\sum_{\pi\in C_3}b_{i(\pi(1))}\otimes b_{i(\pi(2))}\otimes b_{i(\pi(3))},x\rangle.
}
Hence, as $\kfrac13\sum_{\pi\in C_3}b_{i(\pi(1))}\otimes b_{i(\pi(2))}\otimes b_{i(\pi(3))}$ belongs to $W$,
the left-hand side of \rf{14fe15b} is equal to
\dyyz{
\langle
\kfrac13\sum_{\pi\in C_3}b_{i(\pi(1))}\otimes b_{i(\pi(2))}\otimes b_{i(\pi(3))},
\kfrac13\sum_{\rho\in C_3}b_{j(\rho(1))}\otimes b_{j(\rho(2))}\otimes b_{j(\rho(3))}
\rangle,
}
which is equal to the right-hand side of \rf{14fe15b}, as the $b_i$ form an orthonormal basis.

So for any $\varphi:E(G)\to [n]$ and $\psi:E(H)\to [n]$ and
any $u\in V(G)$ and $v\in V(H)$, $\langle \widehat\varphi_u,\widehat\psi_v\rangle$
is equal to 1/3 of the number of bijections $\beta:\delta(u)\to\delta(v)$
such that $\psi\circ\beta=\varphi|_{\delta(u)}$ that preserve the cyclic order.
($\delta(w)$ is the set of edges incident with a vertex $w$.)
This being in conformity with \rf{28se14b}, we have \rf{13ap08c}.
\bx

By the first fundamental theorem of invariant
theory for $O(n)$ (cf.\ [8], and Corollary 2.3 and Lemma 4.5 in [22]),%
\de{4ap07j}{
p_n(\oR[\GG])=\OO(W)^{O(n)},
}
the latter denoting the space of $O(n)$-invariant elements of $\OO(W)$.

\sect{17fe15a}{Proof of the Theorem}

To see necessity in the theorem, let $n\in\oZ_+$ and $(c_{ijk})_{i,j,k=1}^n\in W$ ($=((\oR^n)^{\otimes 3})^{C_3}$).
Then the positive semidefiniteness of $M_{f_c,k}$ follows from
\dyyz{
f_c(G\join{k}H)=
p_n(G\join{k}H)(c)=
\langle d^kp_n(G)(c),d^kp_n(H)(c)\rangle,
}
using Lemma \ref{13ap08d}.

To prove sufficiency, let $f:\GG\to\oR$ be weakly reflection positive.
By Lemma \ref{22fe15a}, $f(\bgcirc)$ belongs to $\oZ_+$.
Set $n:=f(\bgcirc)$.
We show that $f=f_c$ for some $c\in((\oR^n)^{\otimes 3})^{C_3}$.
First:
\dy{20ap08b}{
there is an algebra homomorphism $F:p_n(\oR[\GG])\to\oR$ such that $f=F\circ p_n$.
}
Otherwise, as $p_n$ and $f$ (extended to $\oR[\GG]$) are algebra homomorphisms, there is a $\gamma\in\oR[\GG]$ with $p_n(\gamma)=0$
and $f(\gamma)\neq 0$.
We can assume that $\gamma$ is homogeneous, that is, all graphs
in $\gamma$ have the same number of vertices, $k$ say.
So $\gamma\join{k}\gamma$ has no vertices, that is, it is a polynomial in $\bgcirc$.
As moreover $f(\bgcirc)=n=p_n(\bgcirc)$,
we have $f(\gamma\join{k}\gamma)=p_n(\gamma\join{k}\gamma)=0$,
the latter equality because of Lemma \ref{13ap08d}.
By the weak reflection positivity of $f$ this implies that $f(\gamma\join{k}H)=0$ for each $H\in\GG'$.
Hence, by applying $f$ to both sides of the linearization of \rf{22fe15d} (substituting $\gamma$ for $G$),
$f(\gamma)=0$.
This proves \rf{20ap08b}.

As in [4], \rf{20ap08b} with \rf{4ap07j} implies the existence of $c$ in the complex extension of $W$
satisfying $F(q)=q(c)$ for each $q\in\OO(W)^{O(n)}=p_n(\oR[\GG])$.
To prove that we can take $c$ real, we apply the Procesi-Schwarz theorem [16].
For all $G,H\in\GG$, using Lemma \ref{13ap08d}:
\de{17ap08b}{
F(\langle dp_n(G),dp_n(H)\rangle)
=
F(p_n(G\join{1}H))
=
f(G\join{1}H)
=
(M_{f,1})_{G,H}.
}
Since $M_{f,1}$ is positive semidefinite, \rf{17ap08b} implies
that for each $q\in p_n(\oR[\GG])$:
$F(\langle dq,dq\rangle)\geq 0$, and hence by [16] we can take $c$ real.

Concluding, $f(G)=F(p_n(G))=p_n(G)(c)=f_c(G)$ for each $G\in\GG$, as required.
\bx

\sect{22fe15b}{Uniqueness of $c$}

We finally observe that if $f$ is a partition function,
then $f=f_c$ for some unique $c$, up to the natural action of $O(n)$ on $c$ (which action leaves $f_c$
invariant (cf.\ \rf{4ap07j})).
To see this, let $b\in((\oR^m)^{\otimes 3})^{C_3}$ and
$c\in((\oR^n)^{\otimes 3})^{C_3}$ with $f_b=f_c$.
Then $m=f_b(\bgcirc)=f_c(\bgcirc)=n$.
We show that there exists $U\in O(n)$ such that $b=c^U$
(where $x\mapsto x^U$ is the natural action of $U$ on $x\in W$).

Suppose to the contrary that $b\neq c^U$ for each $U\in O(n)$.
Then the sets
$\{b^U\mid U\in O(n)\}$
and
$\{c^U\mid U\in O(n)\}$
are disjoint compact subsets of $W$.
So, by the Stone-Weierstrass theorem, there exists a polynomial $q\in\OO(W)$
such that $q(b^U)\leq 0$ for each $U\in O(n)$ and $q(c^U)\geq 1$ for each $U\in O(n)$.
As $O(n)$ is compact, we can average $q$ to make it $O(n)$-invariant.
Hence by \rf{4ap07j}, $q\in p_n(\oR[\GG])$, say $q=p_n(\gamma)$ with
$\gamma\in\oR[\GG]$.
Then $f_b(\gamma)=p_n(\gamma)(b)=q(b)\leq 0$ and $f_c(\gamma)=p_n(\gamma)(c)=q(c)\geq 1$.
This contradicts $f_b=f_c$.

\bigskip
\noindent
{\em Acknowledgements.}
We are very thankful to the referees for reading carefully the manuscript and for giving excellent
suggestions to correct the paper and to add helpful explanations.

\section*{References}\label{REF}
{\small
\begin{itemize}{}{
\setlength{\labelwidth}{4mm}
\setlength{\parsep}{0mm}
\setlength{\itemsep}{0mm}
\setlength{\leftmargin}{5mm}
\setlength{\labelsep}{1mm}
}
\item[\mbox{\rm[1]}] D. Bar-Natan, Lie algebras and the four color theorem,
{\em Combinatorica} 17 (1997) 43--52.

\item[\mbox{\rm[2]}] B. Bollob\'as, O. Riordan, 
A polynomial invariant of graphs on orientable surfaces,
{\em Proceedings of the London Mathematical Society} (3) 83 (2001) 513--531.

\item[\mbox{\rm[3]}] S. Chmutov, S. Duzhin, J. Mostovoy, 
{\em Introduction to Vassiliev Knot Invariants},
Cambridge University Press, Cambridge, 2012.

\item[\mbox{\rm[4]}] J. Draisma, D. Gijswijt, L. Lov\'asz, G. Regts, A. Schrijver, 
Characterizing partition functions of the vertex model,
{\em Journal of Algebra} 350 (2012) 197--206.

\item[\mbox{\rm[5]}] S.V. Duzhin, A.I. Kaishev, S.V. Chmutov, 
The algebra of $3$-graphs,
{\em Proceedings of the Steklov Institute of Mathematics} 221 (1998) 157--186.

\item[\mbox{\rm[6]}] M.H. Freedman, L. Lov\'asz, A. Schrijver, 
Reflection positivity, rank connectivity, and homomorphisms of graphs,
{\em Journal of the American Mathematical Society} 20 (2007) 37--51.

\item[\mbox{\rm[7]}] J. Fr\"ohlich, R. Israel, E.H. Lieb, B. Simon, 
Phase transitions and reflection positivity. I. General theory and long range lattice models,
{\em Communications in Mathematical Physics} 62 (1978) 1--34. 

\item[\mbox{\rm[8]}] R. Goodman, N.R. Wallach, 
{\em Symmetry, Representations, and Invariants},
Springer, Dordrecht, 2009.

\item[\mbox{\rm[9]}] J.L. Gross, T.W. Tucker, 
{\em Topological Graph Theory},
Wiley, New York, 1987.

\item[\mbox{\rm[10]}] P. Hanlon, D. Wales, 
On the decomposition of Brauer's centralizer algebras,
{\em Journal of Algebra} 121 (1989) 409--445.

\item[\mbox{\rm[11]}] P. de la Harpe, V.F.R. Jones, 
Graph invariants related to statistical mechanical models:
examples and problems,
{\em Journal of Combinatorial Theory, Series B} 57 (1993) 207--227.

\item[\mbox{\rm[12]}] L. Heffter, 
Ueber das Problem der Nachbargebiete,
{\em Mathematische Annalen} 38 (1891) 477--508.

\item[\mbox{\rm[13]}] M. Kontsevich, 
Vassiliev's knot invariants, in: {\em I.M. Gelfand seminar Part 2},
{\em Advances in Soviet Mathematics} 16 (2) (1993) 137--150.

\item[\mbox{\rm[14]}] R.J. Milgram, R.C. Penner, 
Riemann's moduli space and the symmetric groups,
in: {\em Mapping class groups and moduli spaces of Riemann surfaces},
Contemporary Mathematics 150, American Mathematical Society, Providence, R.I., 1993, pp. 247--290.

\item[\mbox{\rm[15]}] T. Ohtsuki, 
{\em Quantum Invariants ---
A Study of Knots, 3-Manifolds, and Their Sets},
World Scientific, River Edge, N.J., 2002. 

\item[\mbox{\rm[16]}] C. Procesi, G. Schwarz, 
Inequalities defining orbit spaces,
{\em Inventiones Mathematicae} 81 (1985) 539--554.

\item[\mbox{\rm[17]}] N.Yu.\ Reshetikhin, V.G. Turaev, 
Ribbon graphs and their invariants derived from quantum groups,
{\em Communications in Mathematical Physics} 127 (1990) 1--26.

\item[\mbox{\rm[18]}] B.E. Sagan, 
{\em The Symmetric Group: Representations, Combinatorial Algorithms, and Symmetric Functions},
Graduate Texts in Mathematics, Vol. 203, Springer, New York, 2001.

\item[\mbox{\rm[19]}] A. Schrijver, 
Graph invariants in the edge model,
in: {\em Building Bridges --- Between Mathematics and Computer Science}
(M. Gr\"otschel, G.O.H. Katona, eds.), Springer, Berlin, 2008,
pp. 487--498.

\item[\mbox{\rm[20]}] A. Schrijver, 
Graph invariants in the spin model,
{\em Journal of Combinatorial Theory, Series B} 99 (2009) 502--511.  

\item[\mbox{\rm[21]}] A. Schrijver, 
On Lie algebra weight systems for 3-graphs,
{\em Journal of Pure and Applied Algebra} 219 (2015) 4597--4606.

\item[\mbox{\rm[22]}] B. Szegedy, 
Edge coloring models and reflection positivity,
{\em Journal of the American Mathematical Society}
20 (2007) 969--988.

\item[\mbox{\rm[23]}] W.T. Tutte, 
On dichromatic polynominals,
{\em Journal of Combinatorial Theory} 2 (1967) 301--320.

\item[\mbox{\rm[24]}] V.A. Vassiliev, 
Cohomology of knot spaces,
in: {\em Theory of Singularities and its Applications} (V.I. Arnold, ed.),
{\em Advances in Soviet Mathematics} 1 (1990) 23--69.

\end{itemize}
}

\end{document}